\documentclass{amsart}

\usepackage{amsmath,amssymb}

\newtheorem{theorem}{Theorem}
\newcommand{\bt}{\begin{theorem}}
\newcommand{\et}{\end{theorem}}
\newtheorem{lemma}{Lemma}
\newcommand{\bl}{\begin{lemma}}
\newcommand{\el}{\end{lemma}}

\newcommand{\beq}{\begin{equation}}
\newcommand{\eeq}{\end{equation}}
\newcommand{\benum}{\begin{enumerate}}
\newcommand{\eenum}{\end{enumerate}}
\newcommand{\N}{\ensuremath{ \mathbf N }}
\newcommand{\Z}{\ensuremath{\mathbf Z}}

\newcommand{\F}{\ensuremath{\mathbf F }}

\newcommand{\Fq}{\ensuremath{{\mathbf F}_q}}

\newcommand{\mcr}{\ensuremath{ \mathcal R}}

\DeclareMathOperator{\qand}{\quad\text{and}\quad}
\DeclareMathOperator{\qqand}{\qquad\text{and}\qquad}

\DeclareMathOperator{\level}{level}

\title[More differences than multiple sums]{The Haight-Ruzsa method for sets  
with more differences than multiple sums}
\author{Melvyn B. Nathanson}
\address{Department of Mathematics\\Lehman College (CUNY)\\Bronx, NY 10468} \email{melvyn.nathanson@lehman.cuny.edu}

\thanks{Supported in part by a grant from the PSC-CUNY Research Award Program.}

\date{\today}

\begin{document}

\maketitle


\begin{abstract}
Let $h$ be a positive integer and let $\varepsilon > 0$.  The Haight-Ruzsa method produces 
a positive integer $m^*$ and a subset $A$ of the additive abelian group $\Z/m^*\Z$ 
such that the difference set is large in the sense that $A-A = \Z/m^*\Z$ and $h$-fold sumset 
is small in the sense that $|hA| < \varepsilon m^*$.  
This note describes, and in a modest way extends, the Haight-Ruzsa argument, 
and constructs sets with more differences than multiple sums in other additive abelian groups.

\end{abstract}

\section{Sets with more sums than differences}
Let $W$ be an additive abelian group.  
For every subset $A$ of $W$, we define the \emph{difference set} 
\[
A-A = \{ a_1 - a_2: a_1, a_2 \in A\}
\]
and, for every positive integer $h$,  the \emph{$h$-fold sumset} 
\[
hA = \{ a_1+a_2+\cdots + a_h: a_i \in A \text{ for all } i=1,2,\ldots, h \}.
\]
In particular, 
\[
2A = \{ a_1+a_2: a_1, a_2 \in A  \}.
\]
Because 
\[
a_1 + a_2 = a_2 + a_1
\]
but
\[
a_1 - a_2 = - (a_2 - a_1)
\]
it is reasonable to expect that 
\[
|A+A| \leq |A-A|
\]
for ``most'' but not necessarily all finite nonempty subsets $A$ of $W$.  
A set with more sums than differences is called an MSTD set.  
For example, in the additive group \Z\ of integers, the set 
\[
A = \{0,2,3,4,7,11,12,14\}
\]
is an MSTD set: 
\[
A+A =  [0,28]\setminus \{1,20,27\}
\]
\[
A-A = [-14,14]\setminus \{\pm 6, \pm 13\}
\]
and so 
\[
|A+A| = 26 > 25 = |A-A|.
\]
Several families of MSTD sets of integers have been constructed, 
but there is no classification of such sets and many unsolved problems remain
(cf. Hegarty~\cite{hega07a}, Hegarty and Miller~\cite{hega-mill09}, 
Iyer, Lazarov, Miller, and Zhang~\cite{iyer-laza-mill-zhan12,iyer-laza-mill-zhan14}, 
Martin and O'Bryant~\cite{mart-obry07a}, Nathanson~\cite{nath07xyz,nath07yzx}).   

A dual problem is to construct sets with more differences than  multiple sums, 
that is, finite sets $A$ in a group $W$ such that the difference set $|A-A|$ is large 
but the $h$-fold sumset $hA$ is small.  In 1973, Haight~\cite{haig73} proved that for 
all positive integers $h$ and $\ell$ there exists a modulus $m^*$ and a subset $A$ of $\Z/m^*\Z$ 
such that 
\[
A-A = \Z/m^*\Z
\]
but $hA$ omits $\ell$ consecutive congruence classes.  
Recently, Ruzsa~\cite{ruzs16} refined and improved Haight's method, and proved the following:
For every positive integer $h$ and every $\varepsilon > 0$, there exists a modulus $m^*$ 
and a subset $A$ of $\Z/m^*\Z$ 
such that 
\[
A-A = \Z/m^*\Z
\]
and
\[
|hA| < \varepsilon m^*.
\]
The purpose of this note is to describe, and in a modest way extend, the Haight-Ruzsa argument, 
and construct sets with more differences than multiple sums in other additive abelian groups.

\section{Haight-Ruzsa method}
Let $W$ be an abelian group, and let  $f:W \rightarrow W$  be a function, 
not necessarily a homomorphism.  
We define the subset  
\[
A(W,f) = \left\{ w +  f(w) : w \in W  \right\} \cup  \left\{ f(w)) : w \in W  \right\}.
\]
For all $w \in W$, we have 
\[
w = (w+f(w)) - f(w)
\]
and so $A(W,f)$ is a \emph{subtractive basis} for $W$, that is, 
$A(W,f)$ satisfies the difference set identity  
\[
A - A = W.
\]
For every positive integer $h$ and $\varepsilon > 0$, the Haight-Ruzsa method constructs  a 
finite abelian group $W$ and a function $f:W \rightarrow W$ such that 
the $h$-fold sumset of $A(W,f)$ is small in the sense that 
\beq                     \label{HR:ab0}
|hA(W,f)| < \varepsilon |W|.
\eeq

Let $z \in W$.  For every positive integer $h$, 
the element $z$ is in the $h$-fold sumset $hA(W,f)$ if and only if there exist functions  
\beq                     \label{HR:ab1}
\alpha, \beta:W \rightarrow  \{ 0,1,2,\ldots, h \}
\eeq
such that 
\beq                     \label{HR:ab2}
h = \sum_{ w \in W} \left( \alpha(w) + \beta(w) \right)
\eeq
and
\beq                     \label{HR:ab3}
z = \sum_{ w \in W} \left( \alpha(w) (w+f(w)) + \beta(w) f(w) \right).
\eeq
A pair of functions $(\alpha,\beta)$ that satisfies 
conditions~\eqref{HR:ab1} and~\eqref{HR:ab2} is called an \emph{admissible pair for $W$}, 
and the element $z \in hA(W,f)$ defined by~\eqref{HR:ab3} is called the group 
element represented by the admissible pair $(\alpha,\beta)$.  
An element $z \in hA(W,f)$  can be represented by many different admissible pairs.  

An admissible pair $(\alpha,\beta)$ has \emph{level} $\ell$ if
\beq                     \label{HR:ab4}
\ell = \sum_{\substack{w \in W \\ \alpha(w)+\beta(w) \geq 1}} 1.
\eeq
Condition~\eqref{HR:ab2} implies that $\alpha(w)+\beta(w) \geq 1$ for some $w \in W$, 
and so 
\[
\ell \in \{1,2,\ldots, h\}.
\]
If $z$ in $hA$ is represented by the pair $(\alpha,\beta)$, then the level of the pair 
counts the number of $w \in W$ such that at least one of the group elements 
$w+f(w)$ and $f(w)$ appears in the representation~\eqref{HR:ab3} of $z$.

Let $L_{\ell}(W,f)$ be the set of all $z \in W$ such that $z$ can be represented 
by an admissible pair of level at most $\ell$.  We have
\[
L_1(W,f) \subseteq L_2(W,f) \subseteq \cdots \subseteq L_h(W,f) = hA(W,f).
\]
The Haight-Ruzsa method inductively constructs a sequence of groups 
\[
W_1 \subseteq W_2 \subseteq \cdots \subseteq W_h=W
\]
and functions $f_i:W_i \rightarrow W_i$ for $i=1,\ldots, h$ such that 
\[
A(W_1,f_1) \subseteq A(W_2,f_2)\subseteq \cdots \subseteq A(W_h,f_h)
\]
and the sumset $hA(W_h,f_h) = L_h(W_h,f_h)$ satisfies inequality~\eqref{HR:ab0}.

\section{Preliminary}
Let $h$ be a positive integer, and let $\mcr_h$ be the set of commutative rings 
with identity such that, if $R \in \mcr_h$, 
then  $r$ is a unit in $R$ for all $r \in \{1,2,\ldots, h\}$.  
For example, if $p$ is a prime number and $p > h$, then $\Z/p\Z \in \mcr_h$.  
Let $R_0, R_1,\ldots, R_n \in \mcr_h$.  
The direct sum  $\bigoplus_{i=0}^n R_i$ is a ring with identity 
$(1_{R_0}, 1_{R_1},\ldots, 1_{R_n})$, 
and $(r_0, r_1,\ldots, r_n)$ is a unit in $\bigoplus_{i=0}^n R_i$ 
if $r_i \in \{1,\ldots, h\}$ for all $i=0,1,\ldots, n$.  
Thus, $\bigoplus_{i=0}^n R_i \in \mcr_h$.   
If $M_i$ is a finite $R_i$-module for $i=0,1,\ldots, n$, then 
$\bigoplus_{i=0}^n M_i$ is a finite $\bigoplus_{i=0}^n R_i$-module.  
The group $W_h$ will be constructed as the direct sum of a finite number of 
finite $R_i$-modules $M_i$, 
where $R_i \in \mcr_h$ for $i=0,1,\ldots, n$.  

We use the following simple combinatorial inequality.

\bl          \label{HR:lemma:inequality}
Let $M_0,M_1,\ldots, M_n$ be finite sets, 
let 
\[
W = M_0 \times M_1 \times \cdots \times M_n
\] 
and let 
\[
(x_0^*, x_1^*, \ldots, x_n^* ) \in W.
\]
If 
\[
S = \left\{  (x_0,x_1,\ldots, x_n) \in W: x_j = x_j^* \text{ for some } j \in \{0,1,\ldots, n\} \right\}
\]
then
\[
|S| \leq 
|W| \sum_{j=0}^n\frac{1}{|M_j|}.  
\]
If $\varepsilon > 0$ and $|M_j| > (n+1)/\varepsilon$ for all $ j \in \{0,1,\ldots, n\}$, then 
$|S| < \varepsilon |W|$.  
\el

\begin{proof}
The cardinality of $W$ is
\[
|W| = \prod_{i=0}^n |M_i|.
\]
For $j \in \{0,1,\ldots, n\}$, the number of elements $(x_0, x_1, \ldots, x_n) \in W$ with 
$x_j = x_j^*$ is   
\[
\prod_{\substack{ i=0 \\ i \neq j}}^n |M_i| = \frac{|W|}{|M_j|}.
\]
If 
\[
S = \left\{  (x_0,x_1,\ldots, x_n) \in W: x_j = x_j^* \text{ for some } j \in \{0,1,\ldots, n\} \right\}
\]
then
\[
|S| \leq \sum_{j=0}^n\frac{|W|}{|M_j|} = |W| \sum_{j=0}^n\frac{1}{|M_j|}.  
\]
If $\varepsilon > 0$ and $|M_j| > (n+1)/\varepsilon$ for all $ j \in \{0,1,\ldots, n\}$, then 
$|S| < \varepsilon |W|$.  
This completes the proof.  
\end{proof}

\section{Initial step}
Let $h$ be a positive integer and let $\varepsilon > 0$.  
Choose numbers $\varepsilon_1,\ldots, \varepsilon_h$ such that 
\[
0 < \varepsilon_1 < \varepsilon_2 < \cdots < \varepsilon_h < \varepsilon.
\]
For $i= 0,1,2,\ldots, h$, let $R_i \in \mcr_h$, let $M_i$ be a finite $R_i$-module
such that $|M_i|>(h+1)/\varepsilon_1$, and let
\[
W_1 = \bigoplus_{i=0}^h M_i.  
\]
We write the element $w_1 \in W_1$  as an $(h+1)$-tuple  
\beq         \label{HR:w1-tuple}
w_1 = (x_0,x_1,\ldots, x_{h})
\eeq  
where $x_i \in M_i$ for $i= 0,1,\ldots, h$.  
Recall that $h$ is a unit in $R_i$.  For $i = 0, 1, \ldots, h$, we define 
the function $g_i: W_1 \rightarrow M_i$ by 
\[
g_i(w_1) = -\frac{i}{h} x_i
\]
and we define the function $f_1: W_1\rightarrow W_1$ by
\begin{align*}
f_1(w_1) & = (g_0(x_0), g_1(x_1),\ldots, g_h(x_h)) \\ 
& = \left( 0,  -\frac{1}{h} x_1, \ldots, -\frac{i}{h} x_i, \ldots,   -\frac{h-1}{h} x_{h-1}, -x_h \right) 
\end{align*}

Let 
\[
A(W_1,f_1) = \left\{ w_1 + f_1(w_1): w_1 \in W_1 \right\} 
\cup \left\{  f_1(w_1) : w_1 \in W_1 \right\}.
\]
The level 1 set  $L_1(W_1,f_1)$ is the set of all $x \in hA(W_1,f_1)$ of the form 
\[
x = j (w_1 +  f_1(w_1) )+ (h-j)  f_1(w_1)  
= j w_1 + h f_1(w_1)  
\]
for some 
\beq   \label{HR:j-set}
j \in \{0,1,2,\ldots, h\}. 
\eeq
We have
\begin{align*}
 j w_1 & + h f_1(w_1)  \\
& =  j(x_0, x_1,\ldots, x_h) + h \left( 0, -\frac{1}{h} x_1, \ldots, -\frac{i}{h} x_i, \ldots, 
-\frac{h-1}{h} x_{h-1}, -x_h \right)\\
& = (jx_0, jx_1,\ldots, jx_i, \ldots, jx_h) - \left( 0, x_1, \ldots, i x_i, \ldots,  h x_h \right)\\
&  = (jx_0, (j-1)x_1,\ldots,(j-i)x_i, \ldots,  (j-h)x_h).  
\end{align*}
It follows from~\eqref{HR:j-set}  that at least one coordinate 
of this $(h+1)$-tuple is 0.  
Applying Lemma~\ref{HR:lemma:inequality}, we obtain 
\[
|L_1(W_1,f_1)| \leq  |W_1| \sum_{i=0}^h \frac{1}{|M_i|}  
< \varepsilon_1 |W_1|.
\]

\section{Inductive step}
Let $1 \leq k \leq h-1$, and assume that we have a ring $R_0 \in \mcr_h$, a finite  
$R_0$-module $W_k$, and a function $f_k:W_k \rightarrow W_k$ 
such that the sets 
\[
A(W_k,f_k) = \left\{ w_k + f_k(w_k): w_k \in W_k  \right\} 
\cup  \left\{ f_k(w_k): w_k \in W_k  \right\}
\]
and
\[
L_k(W_k, f_k) = \{w_k \in hA_k: \level(w_k) \leq k\}
\]
satisfy 
\[
|L_k(W_k, f_k)| < \varepsilon_k |W_k|.
\]
Because  $W_k$ is a finite set, the number of admissible pairs on $W_k$ is finite.  
Let $n$ be the number of admissible pairs of level exactly $k+1$ 
with respect to $W_k$.  We denote these pairs by $(\alpha_i, \beta_i)$ 
for $i=1,2,\ldots, n$.  It follows from~\eqref{HR:ab2} that 
\[
\alpha_i(w_k) + \beta_i(w_k) \in \{0,1,\ldots, h\}
\]
for all $w_k \in W_k$.

For $i=1,2,\ldots, n$, let $R_i \in \mcr_h$, and 
let $M_i$ be a finite $R_i$-module such that 
\[
|M_i| > \frac{n}{ \varepsilon_{k+1} - \varepsilon_k }.
\]
The set  
\[
W_{k+1} = W_k \oplus  \bigoplus_{i = 1}^{n} M_i 
\]
is a finite module over the  ring 
$R_0 \oplus  \bigoplus_{i = 1}^{n} R_i \in \mcr_h$.
We denote the components of $w_{k+1} \in W_{k+1}$ as follows:  
\[
w_{k+1} = (w_k,x_1,\ldots,  x_{n} ) 
\]
where $w_k \in W_k$ and $x_i \in M_i$ for $i=1,\ldots, n$.  
Define the projection $\pi_0:W_{k+1}\rightarrow W_k$ by 
\[ 
\pi_0(w_{k+1}) = w_k. 
\]
For $i=1,2,\ldots, n$, we define the projection $\pi_i:W_{k+1}\rightarrow M_i$ by 
\[ 
\pi_i(w_{k+1}) = x_i 
\]
and we define the function $g_i: W_{k+1} \rightarrow M_i$ as follows:
If $\pi_0(w_{k+1}) = w_k$  and 
\[
\alpha_i(w_k) + \beta_i(w_k) = 0
\]
then  
\[
g_i(w_{k+1}) = 0. 
\]
If
\[
\alpha_i(w_k) + \beta_i(w_k) \in \{1,\ldots, h\}
\]
then  
\[
g_i(w_{k+1}) = -\frac{\alpha_i(w_k) x_i }{ \alpha_i(w_k) + \beta_i(w_k)  }.  
\]
The function $g_i$ is well defined because $r$ is a unit in $R_i$ for all $r \in \{1,\ldots, h\}$.  

Define the function $f_{k+1}:W_{k+1} \rightarrow W_{k+1}$ 
by
\[
f_{k+1}(w_{k+1}) = \left( f_k(w_k),  g_1(w_{k+1}),\ldots, g_{n}(w_{k+1})  \right).
\]
We have
\begin{align*}
w_{k+1} + f_{k+1}(w_{k+1}) 
& = 
\left( w_k + f(w_k), x_1 + g_1(w_{k+1}), \ldots, x_{n} + g_{n}(w_{k+1}) \right). 
\end{align*}
Consider the set 
\begin{align*}
A&(W_{k+1},f_{k+1}) \\
& = \left\{ w_{k+1}  + f_{k+1} (w_{k+1} ): w_{k+1}  \in W_{k+1} \right\}
\cup  \left\{  f_{k+1} (w_{k+1} ): w_{k+1}  \in W_{k+1} \right\}.  
\end{align*}
Because 
\[
\pi_0\left( f_{k+1}(w_{k+1}) \right) = f(w_k)  
\]
and
\[
\pi_0\left(w_{k+1} +  f_{k+1}(w_{k+1}) \right) = w_k +  f(w_k)
\]
it follows that 
\[
\pi_0(A_{k+1}(W_{k+1},f_{k+1})  ) = A_k(W_k, f_k) .
\]

The $(k+1)$-level set $L_{k+1}(W_{k+1},f_{k+1})$ is the set of all elements 
$z_{k+1}  \in hA_{k+1} $ that can be represented by an admissible pair $(\gamma, \delta)$ 
of level at most $k+1$.  
We define the functions
\[
\hat{\gamma}:W_k \rightarrow \N_0   \qqand \hat{\delta}:W_k \rightarrow \N_0 
\]
as follows:  For $w_k \in W_k$, let 
\[
\hat{\gamma}(w_k) 
= \sum_{\substack{(x_1,\ldots, x_n) \\  \in M_1 \times \cdots \times M_{n}}} 
\gamma(w_k,x_1,\ldots, x_n) 
\]
and 
\[
\hat{\delta}(w_k) 
= \sum_{\substack{(x_1,\ldots, x_n) \\  \in M_1 \times \cdots \times M_{n}}} 
\delta(w_k,x_1,\ldots, x_n) 
\]
We have 
\beq    \label{HR:hat-iden}
\hat{\gamma}(w_k) + \hat{\delta}(w_k) 
= \sum_{\substack{(x_1,\ldots, x_{n}) \\  \in M_1 \times \cdots \times M_{n}}} 
\left( \gamma(w_k,x_1,\ldots, x_{n}) + \delta(w_k,x_1,\ldots, x_{n}) \right) 
\eeq
and 
\begin{align*}
\sum_{w_k \in W_k}  & (\hat{\gamma}(w_k)  + \hat{\delta}(w_k) ) \\
& =  \sum_{w_k \in W_k} 
 \sum_{\substack{(x_1,\ldots, x_{n}) \\ \in M_1 \times \cdots \times M_{n}}} 
 \left( \gamma(w_k,x_1,\ldots, x_{n})  +  \delta(w_k,x_1,\ldots, x_{n}) \right)   \\
& =  \sum_{w_{k+1}\in W_{k+1}} \left( \gamma(w_{k+1}) + \delta(w_{k+1}) \right) \\
& =  \  h.
\end{align*}
Thus, $(\hat{\gamma},\hat{\delta})$ is an admissible pair of functions on $W_k$.

Because $(\gamma, \delta)$ is an admissible pair of functions on $W_{k+1}$ 
of level at most $k+1$, it follows that $\gamma(w_{k+1}) +\delta(w_{k+1}) \geq 1$ 
for at most $k+1$ elements $w_{k+1} \in W_{k+1}$.  
Identity~\eqref{HR:hat-iden} implies that if $w_k \in W_k$ and 
$\hat{\gamma}(w_k) + \hat{\delta}(w_k) \geq 1$, 
then there exists $w_{k+1} \in W_{k+1}$ such that $\pi_0(w_{k+1}) = w_k$ and 
$ \gamma(w_{k+1}) + \delta(w_{k+1})  \geq 1$.  
It follows that $\hat{\gamma}(w_k) +\hat{\delta}(w_k) \geq 1$ 
for at most $k+1$ elements $w_k \in W_k$, and so the pair  $(\hat{\gamma},\hat{\delta})$ 
has level at most $k+1$.  
Similarly, if the pair $(\gamma, \delta)$ has level at most $k$, 
then the pair $(\hat{\gamma},\hat{\delta})$ has level at most $k$.

If $z_{k+1} \in L(W_{k+1}, f_{k+1})$, then $z_{k+1}$ is represented by an admissible pair 
$(\gamma, \delta)$ of level at most $k+1$.  We have  
\begin{align*}
z_{k+1} 
& =   \sum_{w_{k+1} \in W_{k+1}} \left(  \gamma(w_{k+1}) \left(w_{k+1} + f_{k+1}(w_{k+1}) \right) 
  + \delta(w_{k+1}) f_{k+1}(w_{k+1}) \right) \\
& =   \sum_{w_{k+1} \in W_{k+1}} \left( (  \gamma(w_{k+1}) ( w_k + f(w_k), 
x_1 + g_1(w_{k+1}),\ldots, x_n +  g_n(w_{k+1})) \right.  \\
& \hspace{2cm}  \left.   +  \delta(w_{k+1})(  f(w_k),  g_1(w_{k+1}),\ldots,   g_n(w_{k+1}) )   \right).
\end{align*}
For $i=1,\ldots, n$,  we have 
\begin{align*}
\pi_i(z_{k+1}) 
& =   \sum_{w_{k+1} \in W_{k+1}}   
\left(\gamma(w_{k+1}) (x_i + g_i(w_{k+1}) + \delta(w_{k+1})  g_i(w_{k+1}) \right) \\
& =    \sum_{w_{k+1} \in W_{k+1}}   
\left( \gamma(w_{k+1}) x_i + (  \gamma(w_{k+1})  + \delta(w_{k+1}) ) g_i(w_{k+1}) \right)   \in M_i.
\end{align*}
For $i=0$, we have 
\begin{align*}
z_k & =   \pi_0(z_{k+1})\\
 & =   \sum_{w_{k+1} \in W_{k+1}}   
\left( \gamma(w_{k+1}) ( w_k + f(w_k) ) +  \delta(w_{k+1})  f(w_k) \right) \\
 & =   \sum_{w_k \in W_k}  
\sum_{\substack{(x_1,\ldots, x_{n}) \\ \in M_1 \times \cdots \times M_{n}}} 
\left( \gamma(w_k,x_1,\ldots, x_{n}) ( w_k + f(w_k) ) 
+ \delta(w_k,x_1,\ldots, x_{n}) f(w_k)  \right)  \\
& =   \sum_{w_k \in W_k}  \hat{\gamma}(w_k) (w_k +  f(w_k) )+ \hat{\delta}(w_k)  f(w_k) 
\end{align*}
and so $z_k$ is an element of $hA(W_k, f_k)$ that is represented 
by the admissible pair $(\hat{\gamma},\hat{\delta})$ of level at most $k+1$.  
If the admissible pair $(\hat{\gamma},\hat{\delta})$ has level at most $k$, then $z_k \in L_k(W_k,f_k)$.

Suppose that $(\hat{\gamma},\hat{\delta})$ is an admissible pair in $W_k$ of level exactly $k+1$.  
It follows that the admissible pair $(\gamma, \delta)$  in $W_{k+1}$ also has level exactly $k+1$.   
Because $\{ (\alpha_i, \beta_i) : i = 1, \ldots, n \}$ is the set of all admissible pairs 
on $W_k$ of level $k+1$, there is a unique integer $j \in \{1,2,\ldots, n\}$ 
such that $(\hat{\gamma},\hat{\delta}) = (\alpha_j, \beta_j)$.  
Let
\[
S = \{ w_k \in W_k: \alpha_j(w_k) + \beta_j(w_k) \geq 1\}.
\]
There are exactly $k+1$ elements  $w_k \in W_k$ that appear 
in the representation of $z_k$ associated with the admissible pair $(\alpha, \beta)$, 
and so $|S| = k+1$.  
If $w_k \in S$, then   
\begin{align*}
1 & \leq  \alpha_j(w_k) + \beta_j(w_k) = \hat{\gamma}(w_k) + \hat{\delta}(w_k) \\
& = \sum_{\substack{(x_1,\ldots, x_{n}) \\  \in M_1 \times \cdots \times M_{n}}} 
\left( \gamma(w_k,x_1,\ldots, x_{n}) +  \delta(w_k,x_1,\ldots, x_{n}) \right).
\end{align*}
Because the admissible pair $(\gamma, \delta)$  in $W_{k+1}$ has level $k+1$,    
for each $w_k \in S$, there is a unique $n$-tuple 
$(y_1,\ldots, y_n)  \in M_1 \times \cdots \times M_{n}$ such that 
\begin{align*}
\alpha_j(w_k) + \beta_j(w_k) 
& = \hat{\gamma}(w_k) + \hat{\delta}(w_k) \\
& = \gamma(w_k, y_1,\ldots, y_n) +  \delta(w_k, y_1,\ldots, y_n) \\
& \geq 1
\end{align*}
and
\[
 \gamma(w_k, x_1,\ldots, x_n) +  \delta(w_k, x_1,\ldots, x_n) = 0
\]
for all $(x_1,\ldots, x_n) \in M_1 \times \cdots \times M_n$ 
with $(x_1,\ldots, x_n) \neq (y_1,\ldots, y_n)$.
Therefore, 
\begin{align*}
\pi_j(z_{k+1}) 
& =    \sum_{w_{k+1} \in W_{k+1}}   
\left( \gamma(w_{k+1}) x_j + (  \gamma(w_{k+1})  + \delta(w_{k+1}) ) g_j( w_{k+1} ) \right)  \\
& =    \sum_{w_k \in W_k}   \sum_{\substack{ (x_1,\ldots, x_n) \in \\ M_1 \times \cdots \times M_n}} 
 \gamma(w_k, x_1,\ldots, x_n) x_j  \\
& \hspace{0.5cm} +   \sum_{w_k \in W_k}   \sum_{\substack{ (x_1,\ldots, x_n) \in \\ M_1 \times \cdots \times M_n}} 
\left(  \gamma(w_k, x_1,\ldots, x_n)   + \delta(w_k, x_1,\ldots, x_n) \right) g_j(w_{k+1})  \\
& =    \sum_{w_k \in W_k}  \left( \gamma(w_k, y_1,\ldots, y_n) x_j 
+ (  \gamma(w_k, y_1,\ldots, y_n)   + \delta(w_k, y_1,\ldots, y_n) ) g_j(w_{k+1}) \right)  \\
& =    \sum_{w_k \in W_k}  \left( \alpha_j(w_k)  x_j 
+ ( \alpha_j(w_k)  + \beta_j(w_k)  ) g_j(w_{k+1}) \right)  \\
& =    \sum_{w_k \in W_k}  \left( \alpha_j(w_k)  x_j 
+ ( \alpha_j(w_k)  + \beta_j(w_k)  )
 \left(  -\frac{\alpha_j(w_k) x_j }{ \alpha_j(w_k) + \beta_j(w_k)  }\right) \right)  \\
 & = 0.
\end{align*}
To summarize,
if $z_{k+1} \in L_{k+1}(W_{k+1},f_{k+1})$, then $\pi_0(z_{k+1}) \in L_k(W_k,f_k)$ 
or  $\pi_j(z_{k+1}) = 0$ for some $j \in \{1,\ldots, n\}$.
Because $|L_k(W_k,f_k)| < \varepsilon_k |W_k|$, the number of elements  $z_{k+1} \in W_{k+1}$ 
with $\pi_0(z_{k+1})  \in L_k(W_k,f_k)$ is at most 
\[
\varepsilon_k |W_k| |M_1| \cdots |M_n| = \varepsilon_k |W_{k+1}|.
\]
For $j \in \{1,\ldots, n\}$, the number of elements in $W_{k+1}$ with $\pi_j(z_{k+1}) = 0$ is 
\[
|W_k| \prod_{ \substack{i=1\\ i\neq j}}^n |M_i| = \frac{|W_{k+1}|}{|M_j|} 
<  \left( \frac{\varepsilon_{k+1} - \varepsilon_k}{n} \right) |W_{k+1}|
\]
and so 
\[
|L_{k+1}(W_{k+1},f_{k+1}| < \varepsilon_k |W_{k+1}| 
+ \sum_{i=1}^n  \left( \frac{\varepsilon_{k+1} - \varepsilon_k}{n} \right) |W_{k+1}| 
= \varepsilon_{k+1} |W_{k+1}|.
\]
This completes the induction.  With $k=h$, we obtain 
\[
hA(W_h,f_h) = L_h(W_h,f_h) < \varepsilon_h |W_h| < \varepsilon |W_h|.
\]
This completes the proof.

\section{Examples}

Example 1:  
Let $h$ be a positive integer and let $\varepsilon > 0$.  
Let $\mathcal{M}_h$ be the set of positive integers $m$ such that 
every prime divisor of $m$ is greater than $h$. 
If $m \in \mathcal{M}_h$, then every integer $r \in \{1,2,\ldots, h\}$ is a unit in the ring $\Z/m\Z$, 
and so $\Z/m\Z \in \mcr_h$ 
and $\Z/m\Z$ is a finite $\Z/m\Z$-module.  
In the Haight-Ruzsa construction, if we choose modules $M_i = \Z/m_i\Z$ 
with pairwise relatively prime moduli $m_i \in \mathcal{M}_h$, 
then we obtain a finite abelian group  
\[
W_h = \bigoplus_{i\in I} \Z/m_i\Z \cong \Z/m^*\Z
\]
where 
\[
m^* = \prod_{i\in I} m_i 
\]
and a function $f_h:W_h \rightarrow W_h$ such that the set 
\[
A = A(W_h,f_h) = \{ w_h + f_h(w_h):w_h \in W_h\} 
\cup  \{f_h(w_h):w_h \in W_h\}
\]
satisfies $A-A = W_h$ and $|hA| < \varepsilon |W_h|$.  
This is Ruzsa's result, with $m_i$ prime for all $i$ and $m^*$ square-free.  \\

Example 2:
Let $h$ be a positive integer and let $\varepsilon > 0$.  
Let $\F_{q_i}$ be the finite field with $q_i = p_i^{k_i}$ and $p_i > h$.  
In the Haight-Ruzsa construction,  we can choose modules $M_i$ 
that are finite-dimensional vector spaces over the field $\F_{q_i}$ 
of sufficiently large dimension.  
We obtain a finite abelian group $W_h$ that is a direct sum of vector spaces.
If each of these vector spaces is a vector space over the same field $\Fq$, 
then $W_h$ is a vector space over $\Fq$.  \\

Example 3:
Let $h$ be a positive integer and let $\varepsilon > 0$.  
Let $\F_q$ be the finite field with $q = p^k$ and $p > h$, and let $\Fq[t]$ be 
the vector space of polynomials with coefficients in $\Fq$.    
Let $d_0 = 0 < d_1 < d_2 < \cdots  $ be a sufficiently rapidly increasing sequence 
of  integers.  We choose finite subspaces of  $\Fq[t]$ of the form  
\[
\left\{ \sum_{j= d_{i-1}}^{d_i-1} c_jt^j :c_j \in \Fq \right\}.
\]
The Haight-Ruzsa method produces a subspace $W_h$ of the vector space $\Fq[t]$ 
consisting of all polynomials of degree less than $d_h$.

\section{Open problems}
1.  Let $h$ be a positive integer and $\varepsilon > 0$.
Does there exist a prime $p$ and a set $A$ in $\Z/p\Z$ such that 
$A-A = \Z/p\Z$ and $|hA| < \varepsilon p$?
Do such sets exist for all sufficiently large primes $p$?

2.  Let $h$ be a positive integer, let $\varepsilon > 0$, and let $c > 0$.  
Do there exist an additive abelian group $W$ with $|W| > c$ 
and a subset $A$ of $W$ such that $A-A = W$ and the set 
\[
G = \{ 2a_1 + a_2 + \cdots + a_h : a_i \in A \text{ for all } i=1,2,\ldots, h \} 
\]
satisfies
\[
|G|  < \varepsilon |W|?
\]

3.  Let $F(x_1,\ldots, x_h) = \sum_{i=1}^h r_ix_i$ be a linear form with  coefficients 
in a ring $R$, and let $W$ be an $R$-module.  For every subset $A$ of $W$, we define
\[
F(A) = \left\{ \sum_{i=1}^h r_ia_i : a_i \in A \text{ for all } i=1,\ldots, h \right\}.
\]  
It is an open problem to determine the pairs of linear forms $(F,G)$ such that, 
for every $\varepsilon > 0$ and  $c  > 0$, there exists a finite $R$-module $W$ with $|W| > c$, 
$F(A) = W$,  and $|G(A)| < \varepsilon |W|$. 
The Haight-Ruzsa produces modules for which this condition is satisfied 
for the pair of linear forms $(F,G)$, where 
\[
F(x_1,x_2) = x_1 - x_2
\]
and 
\[
G(x_1,\ldots, x_h) = x_1 + \cdots + x_h.
\]

4.  If $F$ and $G$ are polynomials with coefficients in a commutative ring $R$, 
and if $\varepsilon > 0$ and $c > 0$, 
does there exist a finite $R$-algebra $W$ with $|W| > c$ 
and a subset $A$ of $W$ such that $F(A) = W$
and $|G(A)| < \varepsilon |W|$.

\def\cprime{$'$} \def\cprime{$'$} \def\cprime{$'$} \def\cprime{$'$}
\providecommand{\bysame}{\leavevmode\hbox to3em{\hrulefill}\thinspace}
\providecommand{\MR}{\relax\ifhmode\unskip\space\fi MR }
\providecommand{\MRhref}[2]{%
  \href{http://www.ams.org/mathscinet-getitem?mr=#1}{#2}
}
\providecommand{\href}[2]{#2}

\end{document}